\documentclass{article}
\usepackage{latexsym}
\usepackage{amssymb}
\usepackage{hyperref}
\usepackage{amsmath}
\usepackage{lscape}

\usepackage{epsf}
\usepackage{indentfirst}
\usepackage{latexsym}
\usepackage{amssymb}
\usepackage[dvips]{graphicx}
\usepackage{flafter}
\usepackage{caption2}
\usepackage{textcomp}
\usepackage{supertabular}
\usepackage{longtable, booktabs}

\usepackage{multirow}

\usepackage{hhline}
\usepackage{bigstrut,bigdelim}

\usepackage{rotating}
\usepackage{multicol}
\usepackage{multirow}
\usepackage{makeidx} 
\usepackage{epsfig}
\usepackage{fancyhdr}       

\newtheorem{theorem}{Theorem}

\newtheorem{lemma}{Lemma}[section]

\newtheorem{remark}{Remark}

\def\G1{G^\mathcal{C}}

\def\la{{\langle}}

\def\lam{\lambda}

\def\Soc{{\rm Soc}}

\def\la{\lambda}

\def\Soc{{\rm Soc}}

\begin{document}
\title{\bf Classification of flag-transitive $2$-$(v, k, \lambda)$ designs with alternating group $A_n$($n \le 10$) as socle}
\author{\Large Delu Tian\footnote{E-mail: tiandelu@gdei.edu.cn }\ ,
{\Large Qianfen Liao}\footnote{E-mail: liaoqianfen@163.com}\\
 {\small \it  School of Mathematics, Guangdong University of
Education,}\\
{\small\it Guangzhou, Guangdong,  510303,  P. R. China}\\
{\Large Zhilin Zhang\footnote{E-mail: 20241032@gdufe.edu.cn}}\\
 {\small \it  School of Statistics and Mathematics, Guangdong University of }\\
{\small\it   Finance and Economics,  Guangzhou, Guangdong,  510320,  P. R. China}\\
\date{}
}

\maketitle
\begin{abstract}
This paper is devoted to the classification of all flag-transitive point-primitive non-trivial $2$-$(v, k, \lambda)$ designs with
the alternating group  $A_n$($n \le 10$)
 as the  socle  of their automorphism groups, and 87 different designs are obtained up to isomorphism.
 The results of this study further improve the classification theory of designs under the action of almost simple groups,
 and provide reference for the follow-up study of similar problems.

\end{abstract}

MSC2020 Classification: 05B05, 20B25, 20B15, 20B30

Keywords: Flag-transitive;  Point-primitive; Block design; Alternating group

\section{Introduction}

An incidence system $(v, b, r, k, \lambda)$ is called a block design where a set  $\cal P$ of $v$ points is divided into a family $\cal B$ of $b$ distinct subsets (blocks) so that every two points lie in exactly $\lambda$ blocks with  $k$ points in every block, and every point is contained in $r$ blocks. It is also generally required that $k<v$, which is where the {\it ``incomplete"} comes  from in the formal term most often encountered for block designs, {\it balanced incomplete block designs} (BIB designs)\cite{Handbook}.
In a design, a {\it flag} $(\alpha,B)$ is an incident point-block pair such that $\alpha\in B \subseteq \cal B$.
  The design $\cal D$ is called a {\it symmetric design} when $b=v$, otherwise, it is called a {\it nonsymmetric design}.

  The {\it complement} of a design has parameters $(v,b,b-r, v-k, b-2r+\lambda)$.

Despite not being independent the five parameters $v,b, r, k$ and $\lambda$ meet the following three relations:

~~~~~
 $vr=bk$, ~~
 $\lambda(v-1)=r(k-1)$, ~~
  $b\geq v$ (Fisher's inequality).

A BIB design  $\cal D$ is therefore commonly written as  $2$-$(v, k, \lambda)$ design, since $b$ and $r$ are given in terms of $v$, $k$,
 and $\lambda$ by

~~~~~~~~~~~~~~~~~~
$b=\frac{\lambda v(v-1)}{k(k-1)}$,~~~
$r=\frac{\lambda (v-1)}{k-1} $.

 When $2 < k < v-1$ holds, we speak of a {\it non-trivial} 2-design.
 A $2$-$(v, k, \lambda)$ design $\cal D$ is called a {\it  full design}  which consisting of all $k$-subsets.

A group $G$ is {\it almost simple} if it satisfies $T \leq G \leq Aut(T)$ for some simple group $T$ called {\it a socle}.
A permutation $g$ of $\cal P$ that causes a permutation on the blocks is called an {\it automorphism} of $\cal D$.
The symbol $Aut(\cal D)$ represents the {\it full automorphism group} of $\cal D$, which is made up of all automorphisms of $\cal D$.
Each subgroup of $Aut(\cal D)$ is called an {\it automorphism group} of $\cal D$.
The design $\cal D =(\cal P,\cal B)$ is referred to as flag-transitive if $G \leq Aut(\cal D)$ is transitive on the set of flags and point-primitive (or block-primitive) if $G$ is primitive on $\cal P$ (or $\cal B$).  Additional standard notations and definitions are available, for instance, in  \cite{Handbook,AtlasBook,SymmetricDesigns,PermutationGroups}.

In 2013, by using the O'Nan-Scott Theorem, Tian and Zhou \cite{Tian_atmost100} proved that if $D$ is a 2-$(v, k, \lambda)$ symmetric design
with $\lambda \leq  100$ admitting a flag-transitive point-primitive automorphism group $G$, then $G$ must be an almost simple or an affine group.
Subsequently,  all flag-transitive point-primitive symmetric $(v,k,\lambda)$ designs with sporadic socle were fully categorized by them \cite{Tian_sporadic}.
In 2022, Alavi et al.\cite{GCD=1} presented a classification of 2-designs with gcd($r, \lambda$)=1 admitting flag-transitive automorphism groups.
Montinaro  et al.\cite{Liang-Mon,Montinaro1,Montinaro2} have recently classified flag-transitive 2-designs under special $\lambda$.
It is meaningful to consider the classification of designs with almost simple group as the socle.

In the research on the classification of the 2-$(v, k, \lambda)$  designs of the flag-transitive point-primitive group with alternating $A_n$ as socle,
scholars have added various restrictions, mainly focusing on the following situations:

  (i) Symmetric design:

  ~~~ (a) $\lambda\leq 100$(\cite{Dong-vk4,Dong-atmost10,REGUEIRO-Biplanes,Zhou-Dong,Zhu-Tian-Zhou}).

  ~~~ (b) $gcd(k,\lambda)=1$(\cite{Zhu-Guan-Zhou}).

  ~~~ (c) $\lambda \geq gcd(k, \lambda)^2$(\cite{WangZhou2020}).

  (ii) Non-symmetric design:

  ~~~ (a) $\lambda =2$(\cite{Liang2016-nosy2}).

  ~~~ (b) $v< 100$(\cite{{Tang-Chen-Zhang}}).

  ~~~ (c) $gcd(r,\lambda)=1$(\cite{Zhou-Wang2015}).

  ~~~ (d) $\lambda \geq gcd(r, \lambda)^2$(\cite{WangYJ2022}).

  (iii) Symmetry and non-symmetry are considered together:

  ~~~ (a) $\lambda =1$(\cite{BDLNPZ09,Delandtsheer-1}).

  ~~~ (b) $gcd(r-\lambda, k)= 1$(\cite{Zhang2023}).

  ~~~ (c) $\lambda$ is a prime numbe(\cite{Zhang-Chen-Zhou2024}).

  ~~~ (d) $r$ is a prime square $p^2$(\cite{Shen-Chen-Zhou-PrimeSquare}).

All the above work is based on the classification of 2-designs by limiting parameters.
We take a different approach and consider some simple groups as the socle of the automorphism groups
to complete the classification of 2- designs.

 In 2020, Tian \cite{Tian_M11} completely classified flag-transitive point-primitive 2-designs with socle $M_{11}$, and discovered exactly 14 nonisomorphic 2-designs. In 2025, the classification of 2-design  with socle $M_{12}$ was published(\cite{Tian_2025M12}).

 The classification of flag-transitive point-primitive 2-designs with alternating group $A_n$ as socle began in 2019,
 and the case of $Soc(G)=A_5$ was relatively simple\cite{Tian_2020A5}. In 2020, the classification of 2-designs with
 $Soc(G)=A_6$ and $A_7$ was completed,
 and the classification of 2-designs with $Soc(G)=A_n(8\leq n\leq 10)$ was completed one after another.

 In the process of summarizing all the 2-designs of these 6 alternating groups $A_n(5\leq n\leq 10)$ as the socle of the automorphism groups,
 we get the following two conclusions.

\smallskip

\begin{theorem}\label{thm1}
 Let $n \ge 5$, there are only $(n-4)$  flag-transitive point-primitive non-trivial $2$-$(v,k,\lambda)$ designs
  when the symmetric group $S_n$ or the alternating group $A_n$
 acting on $n$  points, and all of which are full designs. Those designs' parameters $(v,b,r, k,\lambda)$ are
$(n,(^n_k),(^{n-1}_{k-1}),k,(^{n-2}_{k-2}))$,
where $k\in \{3,4,\cdots,n-2\}$ and $(^n_k)=\frac{n\cdot (n-1)\cdots (n-k+1)}{k\cdot (k-1)\cdots 1}$.

 \end{theorem}


\begin{theorem}\label{thm2}
 Let $\mathcal{D}$ be a non-trivial  $2$-$(v,k,\lambda)$ design and $G$ be a flag-transitive, point-primitive automorphism group
 of almost simple type, if the socle is $\Soc(G)=A_n$ $(5\leq n \leq 10)$,
 then up to isomorphism there exist 87 designs. one of the following applies:

  (i) When the symmetric group $S_n$ and the alternating group $A_n$ act on a set of $n$  points,
  the designs are full designs, listed in Table Table \ref{tab:v=n}.

 (ii) The designs in other cases, that is, when $\Soc(G)=A_n$ and $G$ acts on $v~ (v \neq n)$, are listed in Table \ref{tab:v<>n}.

 \end{theorem}

\begin{remark}{\rm
\begin{enumerate}
\item[(1)]
  Up to isomorphism, there are 87 different 2-designs, including 84 non-symmetric designs and 3 symmetric designs:
  $(15,7,3)$, $(15,8,4)$ and $(35,18,9)$.

\item[(2)] The parameters of the two designs, which are mutually complementary, are as follows:

$(  6, 3, 2 )$ and itself;~
$(  6, 3, 4 )$ and itself;~
$(  7, 3, 5 )$ and $( 7, 4, 10 )$;~
$(  8, 3, 6 )$ and $( 8, 5, 20 )$;~
$(  8, 4, 15 )$ and itself;~
$(  9, 3, 7 )$ and $( 9, 6, 35 )$;~
$(  9, 4, 21 )$ and $( 9, 5, 35 )$;~
$(  10, 3, 8 )$ and $( 10, 7, 56 )$;~
$(  10, 4, 2 )$ and $( 10, 6, 5 )$;~
$(  10, 4, 4 )$ and $( 10, 6, 10 )$;~
$(  10, 4, 28 )$ and $( 10, 6, 70 )$;~
$(  10, 5, 8 )$ and itself;~
$(  10, 5, 16 )$ and itself;~
$(  10, 5, 56 )$ and itself;~
$(  15, 3, 1 )$ and $( 15, 12, 22 )$;~
$(  15, 5, 4 )$ and $( 15, 10, 18 )$;~
$(  15, 5, 12 )$ and $( 15, 10, 54 )$;~
$(  15, 5, 16 )$ and $( 15, 10, 72 )$;~

$(  15, 6, 10 )$ and $( 15, 9, 24 )$;~
$(  15, 6, 40 )$ and $( 15, 9, 96 )$;~
$(  15, 8, 4 )$ and $( 15, 7, 3 )$.

\item[(3)] There are 21 2-designs are full designs which consisting of all $k$-subsets,
and they are labeled as ${\mathcal D}_1$ to ${\mathcal D}_{21}$ in Table \ref{tab:v=n} and  \ref{tab:v<>n}.

  \end{enumerate} }
\end{remark}

 \begin{table*}[!h]
\begin{center}
\scriptsize
\caption{ Designs with $v=n$$(5\leq n\leq 10)$ and alternating $A_n$  socle}
\label{tab:v=n}
\resizebox{\textwidth}{!}{
\begin{tabular}{cccccccccc}
\hline
{\sc Case} & $G$  & $G_x$  & $v$ & $b$ & $r$ & $k$ & $\lambda$   &  Reference\\
\hline
1 & $A_5$, $S_5$   & $A_4$, $S_4$  & 5 & 10 & 6 & 3 & 3 & ${\mathcal D}_1$ \\
\hline

2 & \multirow{2}*{$A_6$, $S_6$} & \multirow{2}*{$A_5$, $S_5$} & 6 & 20 & 10 & 3  & 4 & ${\mathcal D}_2$ \\
\cline{4-9}

3 &           &         & 6 & 15 & 10 & 4  & 6 & ${\mathcal D}_3$ \\
\hline

4 & \multirow{3}*{$A_7$, $S_7$} & \multirow{3}*{$A_6$, $S_6$}  & 7 & 35 & 15 & 3 & 5  & ${\mathcal D}_4$ \\
\cline{4-9}
5 &           &         & 7 & 21 & 15 & 5  & 10 & ${\mathcal D}_5$ \\
\cline{4-9}
6 &           &         & 7 & 35 & 20 & 4  & 10 & ${\mathcal D}_6$ \\
\hline

7 & \multirow{4}*{$A_8$, $S_8$} & \multirow{4}*{$A_7$, $S_7$}  & 8 & 56 & 21 & 3 & 6  & ${\mathcal D}_7$ \\
\cline{4-9}
8 &           &         & 8 & 28 & 21 & 6  & 15 & ${\mathcal D}_8$ \\
\cline{4-9}
9 &           &         & 8 & 70 & 35 & 4  & 15 & ${\mathcal D}_9$ \\
\cline{4-9}
10 &           &         & 8 & 56 & 35 & 5  & 20 & ${\mathcal D}_{10}$ \\
\hline

11 & \multirow{5}*{$A_9$, $S_9$} & \multirow{5}*{$A_8$, $S_8$}  & 9 & 84 & 28 & 3 & 7  & ${\mathcal D}_{11}$ \\
\cline{4-9}
12 &           &         & 9 & 36 & 28 & 7  & 21 & ${\mathcal D}_{12}$ \\
\cline{4-9}
13 &           &         & 9 & 126 & 56 & 4  & 21 & ${\mathcal D}_{13}$ \\
\cline{4-9}
14 &           &         & 9 & 84 & 56 & 6  & 35 & ${\mathcal D}_{14}$ \\
\cline{4-9}
15 &           &         & 9 & 126 & 70 & 5  & 35 & ${\mathcal D}_{15}$ \\
\hline

16 & \multirow{6}*{$A_{10}$, $S_{10}$} & \multirow{6}*{$A_9$, $S_9$}  & 10 & 120 & 36 & 3  & 8  & ${\mathcal D}_{16}$ \\
\cline{4-9}
17 &           &         & 10 & 45 & 36 & 8  & 28  & ${\mathcal D}_{17}$ \\
\cline{4-9}
18 &           &         & 10 & 210 & 84 & 4  & 28 & ${\mathcal D}_{18}$ \\
\cline{4-9}
19 &           &         & 10 & 120 & 84 & 7  & 56 & ${\mathcal D}_{19}$ \\
\cline{4-9}
20 &           &         & 10 & 252 & 126 & 5  & 56 & ${\mathcal D}_{20}$ \\
\cline{4-9}
21 &           &         & 10 & 210 & 126 & 6  & 70 & ${\mathcal D}_{21}$ \\
\hline

\end{tabular}
}
\end{center}
\end{table*}

 \begin{table*}[htbp]
\begin{center}
\caption{ Designs with $v\neq n$$(5\leq n\leq 10)$ and alternating $A_n$  socle}
\label{tab:v<>n}
\resizebox{\textwidth}{!}{
\begin{tabular}{cccccccccc}
\hline
{\sc Case} & $G$  & $G_x$  & $v$ & $b$ & $r$ & $k$ & $\lambda$   &  Reference\\

\hline
1 & $A_5$   & $D_{10} $  & 6  &  10  &  5  &  3  &  2 & ${\mathcal D}_{22}$ \\
\hline

2 & \multirow{2}*{$S_5$} & \multirow{2}*{$5:4$} & 6  &  20  &  10  &  3  &  4  & ${\mathcal D}_{2}$ \\
\cline{4-9}

3 &           &         & 6  &  15  &  10  &  4  &  6 & ${\mathcal D}_3$ \\
\hline

4 & $S_5$, $A_6$, $S_6$ & $D_{12}$, $F_{36}$, $3^2$:$D_8$   & 10  &  15  &  6  &  4  &  2  & ${\mathcal D}_{23}$ \\
\hline

5 & \multirow{2}*{$A_6$, $S_6$} & \multirow{2}*{$F_{36}$, $3^2$:$D_8$} & 10  &  15  &  9  &  6  &  5  & ${\mathcal D}_{24}$ \\
\cline{4-9}

6 &           &         & 10  &  60  &  18  &  3  &  4 & ${\mathcal D}_{25}$ \\
\hline

7 & $A_6$, $M_{10}$ & $F_{36}$, $3^2$:$Q_8$   & 10  &  36  &  18  &  5  &  8  & ${\mathcal D}_{26}$ \\
\hline

8 & $A_6$, $S_6$, $A_7$, $A_8$ & $S_4$, $S_4\times 2$, $L_2(7)$, $2^3$:$L_3(2)$   & 15  &  15  &  8  &  8  &  4  & ${\mathcal D}_{27}$ \\
\hline

9 & $S_6$, $PGL_2(9)$, $P\Gamma L_2(9)$ & $3^2$:$D_8$, $3^2$:8, $3^2$:[$2^4$]   & 10  &  72  &  36  &  5  &  16  & ${\mathcal D}_{28}$ \\
\hline

10 & \multirow{5}*{$M_{10}$, $PGL_2(9)$, $P\Gamma L_2(9)$} & \multirow{5}*{$3^2$:$Q_8$, $3^2$:8, $3^2$:[$2^4$]}
    & 10  &  30  &  12  &  4  &  4  & ${\mathcal D}_{29}$ \\
\cline{4-9}
11 &           &         & 10  &  30  &  18  &  6  &  10  & ${\mathcal D}_{30}$ \\
\cline{4-9}
12 &           &         & 10  &  120  &  36  &  3  &  8 & ${\mathcal D}_{16}$ \\
\cline{4-9}
13 &           &         & 10  &  45  &  36  &  8  &  28 & ${\mathcal D}_{17}$ \\
\cline{4-9}
14 &           &         & 10  &  180  &  72  &  4  &  24 & ${\mathcal D}_{31}$ \\
\hline

15 & $P\Gamma L_2(9)$ & 10:4   & 36  &  180  &  40  &  8  &  8  & ${\mathcal D}_{32}$ \\
\hline

16 & \multirow{8}*{$A_7$, $A_8$} & \multirow{8}*{$L_2(7)$, $2^3$:$L_3(2)$}
    & 15  &  35  &  7  &  3  &  1  & ${\mathcal D}_{33}$ \\
\cline{4-9}
17 &           &         & 15  &  15  &  7  &  7  &  3  & ${\mathcal D}_{34}$ \\
\cline{4-9}
18 &           &         & 15  &  105  &  28  &  4  &  6 & ${\mathcal D}_{35}$ \\
\cline{4-9}
19 &           &         & 15  &  35  &  28  &  12  &  22 & ${\mathcal D}_{36}$ \\
\cline{4-9}
20 &           &         & 15  &  105  &  42  &  6  &  15 & ${\mathcal D}_{37}$ \\
\cline{4-9}
21 &           &         & 15  &  120  &  56  &  7  &  24 & ${\mathcal D}_{38}$ \\
\cline{4-9}
22 &           &         & 15  &  420  &  84  &  3  &  12 & ${\mathcal D}_{39}$ \\
\cline{4-9}
23 &           &         & 15  &  420  &  168  &  6  &  60 & ${\mathcal D}_{40}$ \\
\hline

24 & \multirow{9}*{$A_{7}$} & \multirow{9}*{$L_2(7)$}
    & 15  &  42  &  14  &  5  &  4  & ${\mathcal D}_{41}$ \\
\cline{4-9}
25 &           &         & 15  &  70  &  28  &  6  &  10  & ${\mathcal D}_{42}$ \\
\cline{4-9}
26 &           &         & 15  &  42  &  28  &  10  &  18 & ${\mathcal D}_{43}$ \\
\cline{4-9}
27 &           &         & 15  &  126  &  42  &  5  &  12 & ${\mathcal D}_{44}$ \\
\cline{4-9}
28 &           &         & 15  &  70  &  42  &  9  &  24 & ${\mathcal D}_{45}$ \\
\cline{4-9}
29 &           &         & 15  &  210  &  56  &  4  &  12 & ${\mathcal D}_{46}$ \\
\cline{4-9}
30 &           &         & 15  &  210  &  84  &  6  &  30 & ${\mathcal D}_{47}$ \\
\cline{4-9}
31 &           &         & 15  &  126  &  84  &  10  &  54 & ${\mathcal D}_{48}$ \\
\cline{4-9}
32 &           &         & 15  &  630  &  168  &  4  &  36 & ${\mathcal D}_{49}$ \\
\hline

33 & \multirow{2}*{$A_7$, $S_7$} & \multirow{2}*{$S_5$, $S_5\times 2$}
    & 21  &  70  &  30  &  9  &  12  & ${\mathcal D}_{50}$ \\
\cline{4-9}
34 &           &         & 21  &  252  &  60  &  5  &  12  & ${\mathcal D}_{51}$ \\
\hline

\multirow{2}*{35} & $A_7$, $S_7$,  & $(A_4\times S_3)$:2, $S_4\times S_3$,
  & \multirow{2}*{35}  &  \multirow{2}*{35}  &  \multirow{2}*{18}
  & \multirow{2}*{18}  &  \multirow{2}*{9}  & \multirow{2}*{${\mathcal D}_{52}$} \\

 &  $A_8$, $S_8$  & $2^4$:$(S_3\times S_3)$, $(S_4\times S_4)$:2
  &   &    &    &   &    &  \\

\hline

36 & \multirow{5}*{$A_8$} & \multirow{5}*{$2^3$:$L_3(2)$}
    & 15  &  168  &  56  &  5  &  16  & ${\mathcal D}_{53}$ \\
\cline{4-9}
37 &           &         & 15  &  280  &  112  &  6  &  40  & ${\mathcal D}_{54}$ \\
\cline{4-9}
38 &           &         & 15  &  168  &  112  &  10  &  72 & ${\mathcal D}_{55}$ \\
\cline{4-9}
39 &           &         & 15  &  280  &  168  &  9  &  96 & ${\mathcal D}_{56}$ \\
\cline{4-9}
40 &           &         & 15  &  840  &  224  &  4  &  48 & ${\mathcal D}_{57}$ \\
\hline

41 & \multirow{2}*{$A_8$} & \multirow{2}*{($A_5\times 3$):2}
    & 56  &  840  &  180  &  12  &  36  & ${\mathcal D}_{58}$ \\
\cline{4-9}
42 &           &         & 56  &  840  &  180  &  12  &  36  & ${\mathcal D}_{59}$ \\
\hline

43 & \multirow{2}*{$S_8$} & \multirow{2}*{$S_5\times S_3$}
    & 56  &  1680  &  360  &  12  &  72  & ${\mathcal D}_{60}$ \\
\cline{4-9}
44 &           &         & 56  &  1680  &  360  &  12  &  72  & ${\mathcal D}_{61}$ \\
\hline

45 & \multirow{4}*{$A_9$, $S_9$} & \multirow{4}*{$S_7$, $S_7\times 2$}
    & 36  &  840  &  140  &  6  &  20  & ${\mathcal D}_{62}$ \\
\cline{4-9}
46 &           &         & 36  &  315  &  140  &  16  &  60  & ${\mathcal D}_{63}$ \\
\cline{4-9}
47 &           &         & 36  &  5040  &  840  &  6  &  120 & ${\mathcal D}_{64}$ \\
\cline{4-9}
48 &           &         & 36  &  5040  &  840  &  6  &  120 & ${\mathcal D}_{65}$ \\
\hline

49 & \multirow{3}*{$A_9$} & \multirow{3}*{$L_2(8)$:$ 3$}
    & 120  &  3360  &  504  &  18  &  72  & ${\mathcal D}_{66}$ \\
\cline{4-9}
50 &           &         & 120  &  10080  &  1512  &  18  &  216  & ${\mathcal D}_{67}$ \\
\cline{4-9}
51 &           &         & 120  &  10080  &  1512  &  18  &  216 & ${\mathcal D}_{68}$ \\
\hline

52 & $S_9$ & $3^3:(2\times S_4)$   & 280  &  11340  &  1296  &  32  &  144  & ${\mathcal D}_{69}$ \\
\hline

53 & \multirow{3}*{$A_{10}$, $S_{10}$} & \multirow{3}*{$S_8$, $S_8\times 2$} &
   45  &  1575  &  420  &  12  &  105 & ${\mathcal D}_{70}$ \\
\cline{4-9}
54 &           &         & 45  &  37800  &  10080  &  12  &  2520  & ${\mathcal D}_{71}$ \\
\cline{4-9}
55 &           &         & 45  &  75600  &  20160  &  12  &  5040 & ${\mathcal D}_{72}$ \\
\hline

56 & $A_{10}$, $S_{10}$  & $(A_7\times 3)$:2, $S_7\times S_3$  & 120  &  33600  &  5040  &  18  &  720  & ${\mathcal D}_{73}$ \\
\hline

57 & $A_{10}$  & $(A_7\times 3)$:2 & 120  &  100800  &  15120  &  18  &  2160  & ${\mathcal D}_{74}$ \\
\hline

58 & \multirow{11}*{$A_{10}$, $S_{10}$} & \multirow{11}*{$(A_5\times A_5)$:4, $(S_5\times S_5)$:2} &
   126  &  4725  &  225  &  6  &  9   & ${\mathcal D}_{75}$ \\
\cline{4-9}
59 &           &         & 126  &  2100  &  600  &  36  &  168  & ${\mathcal D}_{76}$ \\
\cline{4-9}
60 &           &         & 126  &  18900  &  900  &  6  &  36  & ${\mathcal D}_{77}$ \\
\cline{4-9}
61 &           &         & 126  &  37800  &  1800  &  6  &  72  & ${\mathcal D}_{78}$ \\
\cline{4-9}
62 &           &         & 126  &  14175  &  1800  &  16  &  216  & ${\mathcal D}_{79}$ \\
\cline{4-9}
63 &           &         & 126  &  75600  &  3600  &  6  &  144  & ${\mathcal D}_{80}$ \\
\cline{4-9}
64 &           &         & 126  &  151200  &  7200  &  6  &  288  & ${\mathcal D}_{81}$ \\
\cline{4-9}
65 &           &         & 126  &  56700  &  7200  &  16  &  864  & ${\mathcal D}_{82}$ \\
\cline{4-9}
66 &           &         & 126  &  25200  &  7200  &  36  &  2016  & ${\mathcal D}_{83}$ \\
\cline{4-9}
67 &           &         & 126  &  25200  &  7200  &  36  &  2016  & ${\mathcal D}_{84}$ \\
\cline{4-9}
68 &           &         & 126  &  113400  &  14400  &  16  &  1728  & ${\mathcal D}_{85}$ \\
\cline{4-9}
\hline

69 & $S_{10}$  & $S_7\times S_3$ & 120  &  201600  &  30240  &  18  &  4320  & ${\mathcal D}_{86}$ \\
\hline

70 & $S_{10}$  & $(S_5\times S_5)$:2 & 126  &  604800  &  28800  &  6  &  1152  & ${\mathcal D}_{87}$ \\
\hline

\end{tabular}
}
\end{center}
\end{table*}

\section{Some Preliminary Results}

We present some preliminary results in this section that are used throughout this paper.

\begin{lemma}  \label{FPB}  \quad
Let $\mathcal{D}=(\cal P, \cal B)$ be a non-trivial  $2$-design and $G \leq Aut(\cal D)$.
 The following three claims are equal for any point $x\in \cal P$ and block $B \in \cal B$:

 (i) $G$ acts flag-transitively on $\cal D$;

 (ii)   $G$ acts point-transitively on $\cal D$ and $G_x$ acts transitively on $B(x)$,
  where $B(x)$ denotes the set of all blocks which are incident with  $x$;

 (iii) $G$ acts block-transitively on $\cal D$ and $G_B$ acts transitively on the points of $B$.

\end{lemma}

 \begin{lemma}  \label{parameter}  \quad
 Let $\mathcal{D}=(\cal P, \cal B)$ be a non-trivial flag-transitive  $2$-$(v,k,\lambda)$  design and $G \leq Aut(\cal D)$.
 Then the following hold:

  (i) $r>\lambda$, $r \geq k$, $r^2 > \la v$;

  (ii) $b\ |\ |G|$, $r\ |\ |G_x|$, where $G_x$ is any point-stabiliser of $G$.

\end{lemma}

\textbf{Proof.} $(i)$  $\mathcal{D}$ is a non-trivial $2$-$(v,k,\lambda)$ design, hence $2 < k < v-1$.
 From equation $r=\frac{\lambda (v-1)}{k-1} $, we get  $r>\la$.
 Fisher's inequality $b\geq v$ implies that $r\geq k$, then
 $$r^2\geq rk >rk+(\la -r)=r(k-1)+\la =\la(v-1)+\la= \la v.$$

$(ii)$ According to Lemma \ref{FPB}, $G_B$ acts transitively on $B$, and $G_x$ acts transitively on $B(x)$ ,
 so $b\ |\ |G|$  and $r\ |\ |G_x|$ holds.
$\hfill\square$

\newpage

 \section{Proof of Theorem 1}

 Both symmetric groups and alternating groups exhibit strong transitivity properties.

  \begin{lemma}(\cite{PermutationGroups}) \label{PermutationGroups}  \quad
 The symmetric group $S_n$ and the alternating group $A_n$ are respectively $n$-transitive and $(n-2)$-transitive
 on $\cal P = $ $\{ 1, 2, \dots, n\}$.
 \end{lemma}

In fact, $A_n$ and $S_n$ are the only subgroups of $S_n$ that are $(n-2)$-transitive on $\cal P = $ $\{ 1, 2, \dots, n\}$.

\begin{lemma}\label{n-2_Transitive}  \quad
Let $n \geq 5$ and let $G$ be a subgroup of $S_n$ that is $(n-2)$-transitive on $\cal P = $ $\{ 1, 2, \dots, n\}$. Then either $G = A_n$ or $G = S_n$.
\end{lemma}

\textbf{Proof.}
Since $G$ is $(n-2)$-transitive, its action on the set of ordered $(n-2)$-tuples of distinct elements from $\cal P$ is transitive. In particular, the size of the orbit of any such tuple is equal to the number of ordered $(n-2)$-tuples, which is
\[
n(n-1)\cdots 4\cdot 3 = \frac{n!}{2}.
\]
Now, fix the ordered tuple $T = (1, 2, \dots, n-2)$. Let $H = G_T$ be the stabiliser of this tuple in $G$. By the orbit-stabiliser theorem, we have
\[
|G| = \frac{n!}{2} \cdot |H|.
\]
Since $G$ is a subgroup of $S_n$, we have $|G| \leq n!$. Therefore,
\[
\frac{n!}{2} \cdot |H| \leq n!.
\]
Hence, $|H|$ is either 1 or 2.

If $|H| = 2$. Then $|G| = n!$, so $G = S_n$.

If $|H| = 1$. Then $|G| = \frac{n!}{2}$. Thus, $G$ is a subgroup of $S_n$ of index 2.
 For $n \geq 5$, since the commutator subgroup of $S_n$ is $A_n$ and any index 2 subgroup is normal and contains the commutator subgroup,
  the only subgroup of $S_n$ of index 2 is $A_n$. Therefore, $G = A_n$.
$\hfill\square$

\begin{lemma}\label{FullDesign}  \quad
Let $n \geq 5$ and  $G$ be an $(n-2)$-transitive permutation group on $\cal P = $ $\{ 1, 2, \dots, n\}$.
If $G$ preserves a $2$-$(n,k,\lambda)$ design $\cal D =(\cal P,\cal B)$,
then $\mathcal{D}$ is the full design which consisting of all $k$-subsets of $\cal P$.
\end{lemma}

\textbf{Proof.}
Since $G$ is $(n-2)$-transitive on $\cal P$, by Lemma \ref{n-2_Transitive},
$G$ must be either the alternating group $A_n$ or the symmetric group $S_n$.

Now, $G$ acts transitively on the set of all $k$-subsets of $\cal P$ for any $k$ with $3 \leq k \leq n-2$.

Since $G$ preserves the design $D$, the set of blocks $\mathcal{B}$ must be a union of orbits of $G$ acting on the set of $k$-subsets.
But since $G$ acts transitively on all $k$-subsets, the only possible orbits are the empty set and the entire set of $k$-subsets.
Since $D$ is a $2$-design, it must contain at least one block, so $\mathcal{B}$ cannot be empty.
Therefore, $\mathcal{B}$ must be the set of all $k$-subsets of $\cal P$, i.e., the design $D$ is the full design.
$\hfill\square$

The existence of a $2$-$(n,k,\lambda)$ design for any $k$ in the given range is guaranteed by taking the full design.

\begin{lemma}\label{Existence}  \quad
Let $n \geq 5$ and  $G$ be an $(n-2)$-transitive permutation group on $\cal P =$ $ \{1, 2, \dots, n\}$.
Then for any integer $k$ with $3 \leq k \leq n-2$, there exists a $2$-$(n,k,\lambda)$ design $\cal D =$ $(\cal P, \mathcal{B})$ preserved by $G$.
\end{lemma}

\textbf{Proof.}
By Lemma \ref{n-2_Transitive} and Lemma \ref{FullDesign},  $G=A_n$ or $S_n$ acts on $\cal P =$ $ \{1, 2, \dots, n\}$,
and if there is a design, it must be a full design.

For any $k\in \{3,4,\cdots, n-2\}$, consider the set $\mathcal{B}$ of all $k$-subsets of $\cal P$. We first show that $(\cal P, \mathcal{B})$ is a $2$-$(n,k,\lambda)$ design with $\lambda = \binom{n-2}{k-2}$. For any two distinct points $x, y \in \cal P$, the number of $k$-subsets containing both $x$ and $y$ is $\binom{n-2}{k-2}$, since we choose the remaining $k-2$ points from the $n-2$ points other than $x$ and $y$. Thus, the design condition is satisfied.

Next, we show that $G$ preserves $\mathcal{B}$. Since $G$ is either $A_n$ or $S_n$, it acts transitively on the set of all $k$-subsets of $\Omega$. That is, for any two $k$-subsets $B_1$ and $B_2$, there exists $g \in G$ such that $B_1^g = B_2$. In particular, for any $B \in \mathcal{B}$ and any $g \in G$, we have $B^g \in \mathcal{B}$ because $B^g $ is also a $k$-subset. Therefore, $G$ preserves the design.

Hence, for any $k$ with $3 \leq k \leq n-2$, the full design is a $2$-$(n,k,\lambda)$ design preserved by $G$.
$\hfill\square$

\begin{lemma}\label{flag-Tr}
Let $n \geq 5$ and  $G$ be an $(n-2)$-transitive permutation group on $\cal P = $ $\{ 1, 2, \dots, n\}$.
Then the group $G$ acts transitively on the set of flags $\mathcal{F}$=$\{(\alpha, B)\ |\ \alpha \in B, B\subseteq \cal B\}$.
\end{lemma}

\textbf{Proof.}
 Consider the stabiliser $G_x$ of $x$ in $G$.
By Lemma \ref{n-2_Transitive}, $G_x$ is isomorphic to either $A_{n-1}$ or $S_{n-1}$, and in either case, it acts transitively on $B(x)$, the set of blocks containing $x$.
This is because $B(x)$  is in one-to-one correspondence with the set of $(k-1)$-subsets of
 $\cal P$$ \setminus \{x\}$, and $G_x$ acts transitively on that set (since $A_{n-1}$ and $S_{n-1}$ are transitive on
  $(k-1)$-subsets for $3 \leq k \leq n-2$).

 Since $G$ acts point-transitively on $\cal D$ and $G_x$ acts transitively on $B(x)$, by Lemma \ref{FPB},
 the group $G$ acts flag-transitively on $\cal D$.
$\hfill\square$

\begin{lemma}\label{Point-Pr}
Let $n \geq 5$ and  $G$ be an $(n-2)$-transitive permutation group on $\cal P = $ $\{ 1, 2, \dots, n\}$.
Then the group $G$ acts primitively on $\cal P = $ $\{ 1, 2, \dots, n\}$.
\end{lemma}

\textbf{Proof.}
 For $n \geq 5$,  $(n-2)$-transitive implies $2$-transitive.
 By Theorem 9.6(\cite{PermutationGroups}), every $2$-transitive group is primitive,
 thus the group $G$ acts primitively on $\cal P = $ $\{ 1, 2, \dots, n\}$.
$\hfill\square$

\textbf{Proof of Theory \ref{thm1}.}
By Lemma \ref{PermutationGroups} to \ref{Point-Pr}, for $n \geq 5$ and $3\leq k \leq n-2$,
 there are only $(n-4)$  flag-transitive point-primitive non-trivial $2$-$(v,k,\lambda)$ designs
  when the symmetric group $S_n$ or the alternating group $A_n$
 acting on $n$  points, and all of which are full designs.

In the process of proving Lemma \ref{Existence}, $\lambda = \binom{n-2}{k-2}$ has been obtained.
Substituting $v=n$ and $\lambda = \binom{n-2}{k-2}$ into the formula
$b=\frac{\lambda v(v-1)}{k(k-1)}$,
$r=\frac{\lambda (v-1)}{k-1} $,
we can get $b = \binom{n}{k}$,  $r = \binom{n-1}{k-1}$.
$\hfill\square$

 \section{Proof of Theorem 2}

 In the next 2 subsections, we prove Theorem \ref{thm2}.

 \subsection{Getting Possible Parameters of 2-Designs}
 In this subsection, we will get all possible parameters of  2-$(v,k,\lambda)$  designs.

  \begin{lemma} \label{divisible}  \quad
 Let $G_x$ be a maximal subgroup of $G$.
 If $r\ |\ |G_x|$, then $b\ |\ |G|$, but not vice versa.
 \end{lemma}

\textbf{Proof.}
Firstly, from $r\ |\ |G_x|$, it can be inferred that $b\ |\ |G|$ holds.
 From $[G:G_x]=v$, we can get $v=\frac{|G|}{|G_x|}$.
 Substitute it into $bk=vr$.
 By the hypothesis, $r\ |\ |G_x|$, it follows that $bk(\frac{|G_x|}{r})=|G|$,
 and hence $b\ |\ |G|$ holds.

 Secondly,  if $b\ |\ |G|$, it does not necessarily follow that $r\ |\ |G_x|$.
 We can give a counterexample.
 Let $G=A_6$, $G_x\cong F_{36}$, then $|G|=360$, $|G_x|=36$ and $v=10$.
Assuming  $k = 4$, $\lambda = 24$, we can calculate $b=180$, $r=72$.
 Obviously, $b\ |\ |G|$ holds, but $r\ |\ |G_x|$ does not.
$\hfill\square$

\begin{remark}
{\rm
\begin{enumerate}
\item[(1)]
This lemma does not conflict with Lemma \ref{parameter}(ii), and $b\ |\ |G|$ and $r\ |\ |G_x|$ may not be established at the same time
without knowing whether the design $\cal D$ is flag-transitive.
According to Lemma \ref{divisible}, when looking for possible design parameters,
we assume that $r\ |\ |G_x|$ holds, so it is unnecessary to verify $b\ |\ |G|$.
\item[(2)]
The counterexample in the proof of Lemma \ref{divisible} also demonstrates that for
$G=A_6$, a design with parameters (10,180,72,4,24) cannot be constructed. However, for
$G=M_{10}$, such a design is constructible, as evidenced by design
${\mathcal D}_{31}$ in Table \ref{tab:v<>n}.
 \end{enumerate}
 }
 \end{remark}

  \begin{lemma}(\cite{Wilson}) \label{Wilson}  \quad
 For $n=5$ or $n \geq 7$, the automorphism group $Aut(A_n)$ $\cong S_n$, while $Aut(A_6) \cong P\Gamma L_2(9)$.
 \end{lemma}

 Assume that there is a non-trivial  2-design $\cal D$ admitting a flag-transitive and point-primitive almost simple automorphism group
 $G$ with socle $A_n$.
 According to lemma \ref{Wilson}, when $n=6$, $G=A_6,S_6,M_{10},PSL_2(9),P\Gamma L_2(9)$,
and when $n=5,7,8,9,10 $,   $G=A_n,S_n$.

  \begin{lemma}(\cite{PermutationGroups}) \label{Max-Pri}  \quad
 Let $x\in \Omega$, $|\Omega|> 1$. A transitive group $G$ on $\Omega$ is primitive
 if and only if $G_x$ is a maximal subgroup of $G$.
 \end{lemma}

 If $M$ is any maximal subgroup of $G$, then the permutation action of $G$ on the cosets of $M$ is primitive,
 so $G$ embeds as a primitive subgroups of $S_m$, where $m=[G:M]$.

 According to Lemma \ref{Max-Pri}, if and only if the stabiliser $G_x$ is a maximum subgroup of $G$, where $x\in\cal P$,
 then $G$ is point-primitive on $\cal P$.
 Consequently, $v=[G:G_x]$.
 In the ATLAS, the maximal subgroups of $G$ are listed \cite{AtlasBook}.

We calculate all possible parameters $(v,b,r, k,\lam)$ that meet the requirements listed below:

 \begin{enumerate}
    \item[(i)] \,
  $G\in \{A_n,S_n,M_{10},PSL_2(9),P\Gamma L_2(9)\}$ with $5\leq n\leq 10$,
  and $G_x$ is one of its maximal subgroups;
    \item[(ii)] \,
  $v=[G:G_x]$;
      \item[(iii)] \,
  $2<k<v-1$;
 \item[(iv)] \,
 $r\ |\ |G_x|$,  $r > \lambda$  and $r^2>\lambda v$ (Lemma \ref{parameter});
 \item[(v)] \,
  $bk=vr$,~    $\lambda(v-1)=r(k-1)$;
   \item[(vi)] \,
  $v\leq b \leq \binom{v}{k}$.

 \end{enumerate}

We obtain 2091 5-tuples of parameters $(v,b,r, k,\lambda)$ with the help of the computer algebra system {\sf GAP} \cite{GAP}.
The numbers of possible 5-tuples  corresponding to group $G$ are listed in Table \ref{tab:Number}.

\begin{table}[htbp]
\centering
\large
\caption{The numbers of possible 5-tuples $(v,b,r, k,\lambda)$}
\label{tab:Number}
\begin{tabular}{lr|lr|lr}
\hline
$G$       & Sum & $G$               & Sum & $G$       & Sum \\
\hline
$A_5$     & 5   & $\text{PSL}_2(9)$ & 22  & $S_8$     & 100 \\
$S_5$     & 6   & $\text{P}\Gamma\text{L}_2(9)$ & 26  & $A_9$     & 395 \\
$A_6$     & 19  & $A_7$             & 120 & $S_9$     & 447 \\
$S_6$     & 25  & $S_7$             & 101 & $A_{10}$  & 297 \\
$M_{10}$  & 22  & $A_8$             & 157 & $S_{10}$  & 349 \\
\hline
\end{tabular}
\end{table}

These possible 5-tuples parameters are verified one by one, and most of them are eliminated in the following three steps.

 \begin{enumerate}
    \item[Step(i)] \,
  According to Lemma \ref{FPB}, $G$ acts on $\cal D$ in a block-transitive manner.
 Consequently, the subgroup $G_B$ of $G$ has the index $b=[G:G_B]$.
It is simple to determine whether there is at least one subgroup with index $b$ of $G$
 by using the {\sc Magma}-command
 {\tt Subgroups(G:OrderEqual:=n)} where $n=|G|/b$ \cite{Magma}.
    \item[Step(ii)] \,
  There is at least one orbit $O$ of $G_B$ with size $k$ and $|O^G|=b$
 because $G_B$ acts transitively on the points of $B$. Check whether there is such an orbit.
     \item[Step(iii)] \,
 For any two points, they must coexist in $\lambda$ different blocks. Check whether this value is a fixed value.
 \end{enumerate}

We take two cases where $S_9$ acts on 280 points as examples.

\textbf{Ex.1}
$(v,b,r, k,\lambda)$=$( 280, 2880, 648, 63, 144 )$.

There is no subgroup with index $2880$ of $G$, so this parameters can be eliminated by Step(i).
$\hfill\square$

\textbf{Ex.2}
$(v,b,r, k,\lambda)$=$( 280, 11340, 1296, 32, 144 )$.

The symmetric group $S_9$ contains 17 conjugacy classes of  subgroups
  with  index 11340,  and analyze them one by one.

 We list the serial number of the conjugate class of subgroups in the first column in Table \ref{tab:Ex}.
 and the orbital lengths under the action of the conjugate class in the second column.
 If the orbital length is 32, then the number of elements in the set generated by this orbital under the action of group $G$ is indicated in parentheses.
The notation $s^t$ indicates that the degree $s$ appears with multiplicity $t$.
In the third column, the treatment method for this situation is given.
$\hfill\square$

 \begin{table*}[htbp]
\begin{center}
\caption{ Analysis of possible 5-tuples $( 280, 11340, 1296, 32, 144 )$ }
\label{tab:Ex}
\resizebox{\textwidth}{!}{
\begin{tabular}{ccc}
\hline
Conjugate class     &  Orbital lengths    &  Reference         \\
\hline
 1    &  $8^9$, $16^7$, $32(2835)^3$             &     Step(ii)   \\

 2    &  $4^6$, $8^{10}$, $16^5$, $32(5670)^3$    &     Step(ii)  \\

 3    &  $8^{9}$, $16^7$, $32(5670)^3$           &     Step(ii)  \\

 4    &  $8^{5}$, $16^5$, $32(2835)$, $32(5670)^4$    &     Step(ii)  \\

 5    &  $4^{10}$, $8^{6}$, $16^8$, $32(2835)^2$    &     Step(ii)  \\

 6    &  $4^2$, $8^{10}$, $16^4$, $32(5670)^4$    &     Step(ii)  \\

 7    &  $8^{9}$, $16^3$, $32(2835)^5$    &     Step(ii)  \\

 8    &  $8^{5}$, $16^5$, $32(5670)^5$    &     Step(ii)  \\

 9    &  $4^2$, $8^{6}$, $16^6$, $32(5670)^4$    &     Step(ii)  \\

 10   &  $8$, $16^7$, $32(2835)^2$, $32(5670)$, $32(11340)^2$    &  Step(ii),   Step(iii)  \\

 11   &  $2^2$, $4^7$, $8^{9}$, $16^7$, $32(11340)^2$    &     Step(iii)  \\

 12   &  $4^4$, $8^{7}$, $16^7$, $32(5670)$, $32(11340)^2$   &    Step(ii),   Step(iii)  \\

 13   &  $8^{3}$, $16^8$, $32(2835)$, $32(5670)$, $32(11340)^2$   &     Step(ii),${\mathcal D}_{69}$  \\

 14   &  $4^2$, $8^{6}$, $16^6$, $32(5670)^2$, $32(11340)^2$   &     Step(ii),   Step(iii)\\

 15   &  $4^2$, $8^{6}$, $16^8$, $32(5670)$, $32(11340)^2$     &     Step(ii),   Step(iii)  \\

 16   &  $8^{3}$, $16^8$, $32(2835)$, $32(5670)$, $32(11340)^2$     &  Step(ii),   Step(iii)  \\

 17   &  $4^4$, $8^{5}$, $16^8$, $32(5670)$, $32(11340)^2$    &     Step(ii),   Step(iii)  \\

\hline

\end{tabular}
}
\end{center}
\end{table*}

 \subsection{Basic blocks of 2-Designs }

     After talking about the ``knockout" in the last subsection,
     the remaining parameters can get a total of 87 designs
     in Table \ref{tab:v=n} and Table \ref{tab:v<>n} up to isomorphism.

     In fact, when $n=5,6,7,8,9,10$, the 21 designs in Table \ref{tab:v=n} can also be obtained directly from Theorem \ref{thm1}.
    Because they are all full designs, each $k$-tuple can be used as a basic block of $2$-$(n,k, \lambda)$, so it will not be listed here.

   We list the basic blocks of the non-full designs list in Table \ref{tab:v<>n}.
   In the sense of isomorphism, only the basic block of the design under the action of the group $G$ in the first place are listed.

 \begin{table*}[htbp]
\begin{center}
\caption{ Up to isomorphism, the basic block of each design}
\label{tab:Basic block}
\resizebox{\textwidth}{!}{
\begin{tabular}{ccccc}
\hline
{\sc No.} & $G$  & $G_x$  & Basic block    &  Design  \\

\hline
1 & $A_5$   & $D_{10}$  &   \{ 1, 4, 6 \}   & ${\mathcal D}_{22}$ \\
\hline

2 & $S_5$   & $D_{12}$  &  \{ 5, 7, 8, 10 \}   & ${\mathcal D}_{23}$ \\
\hline

3 & \multirow{3}*{$A_6$} & \multirow{3}*{$F_{36}$}  &    \{ 2, 3, 5, 6, 7, 10 \}     & ${\mathcal D}_{24}$ \\
\cline{4-5}

4 &           &         &    \{ 2, 3, 4 \}     & ${\mathcal D}_{25}$ \\
\cline{4-5}

5 &           &         &    \{ 1, 2, 5, 8, 9 \}     & ${\mathcal D}_{26}$ \\
\hline

6 & $A_6$     & $S_4$    &    \{ 4, 5, 6, 7, 12, 13, 14, 15 \}     & ${\mathcal D}_{27}$ \\
\hline

7 & $S_6$     & $3^2$:$D_8$     &     \{ 1, 4, 5, 6, 9 \}    & ${\mathcal D}_{28}$ \\
\hline

8 & \multirow{3}*{$M_{10}$} & \multirow{3}*{$3^2$:$Q_8$}
    &     \{ 3, 5, 7, 10 \}      & ${\mathcal D}_{29}$ \\
\cline{4-5}
9 &           &         &     \{ 2, 4, 5, 6, 8, 10 \}           & ${\mathcal D}_{30}$ \\
\cline{4-5}
10 &           &         &        \{ 1, 2, 5, 9 \}          & ${\mathcal D}_{31}$ \\
\hline

11 & $P\Gamma L_2(9)$ & 10:4   &   \{ 4, 5, 17, 22, 27, 31, 32, 35 \}  & ${\mathcal D}_{32}$ \\
\hline

12 & \multirow{17}*{$A_7$} & \multirow{17}*{$L_2(7)$}
    &    \{ 3, 12, 15 \}         & ${\mathcal D}_{33}$ \\
\cline{4-5}
13 &           &         &      \{ 1, 2, 3, 8, 9, 10, 11 \}              & ${\mathcal D}_{34}$ \\
\cline{4-5}
14 &           &         &     \{ 2, 3, 8, 9 \}       & ${\mathcal D}_{35}$ \\
\cline{4-5}
15 &           &         &       \{ 1, 3, 4, 5, 6, 7, 9, 11, 12, 13, 14, 15 \}             & ${\mathcal D}_{36}$ \\
\cline{4-5}
16 &           &         &        \{ 2, 5, 9, 11, 12, 14 \}            & ${\mathcal D}_{37}$ \\
\cline{4-5}
17 &           &         &       \{ 4, 6, 7, 12, 13, 14, 15 \}              & ${\mathcal D}_{38}$ \\
\cline{4-5}
18 &           &         &        \{ 4, 6, 11 \}              & ${\mathcal D}_{39}$ \\
\cline{4-5}
19 &           &         &        \{ 2, 3, 6, 11, 14, 15 \}              & ${\mathcal D}_{40}$ \\
\cline{4-5}
20 &           &         &       \{ 4, 11, 12, 13, 14 \}                & ${\mathcal D}_{41}$ \\
\cline{4-5}
21 &           &         &      \{ 2, 3, 9, 10, 13, 15 \}                 & ${\mathcal D}_{42}$ \\
\cline{4-5}
22 &           &         &        \{ 1, 2, 4, 6, 10, 11, 12, 13, 14, 15 \}            & ${\mathcal D}_{43}$ \\
\cline{4-5}
23 &           &         &       \{ 2, 4, 5, 9, 10 \}                   & ${\mathcal D}_{44}$ \\
\cline{4-5}
24 &           &         &      \{ 1, 2, 3, 4, 5, 9, 11, 13, 14 \}                & ${\mathcal D}_{45}$ \\
\cline{4-5}
25 &           &         &       \{ 4, 7, 14, 15 \}               & ${\mathcal D}_{46}$ \\
\cline{4-5}
26 &           &         &       \{ 1, 6, 11, 13, 14, 15 \}               & ${\mathcal D}_{47}$ \\
\cline{4-5}
27 &           &         &         \{ 1, 2, 4, 5, 6, 7, 8, 10, 14, 15 \}           & ${\mathcal D}_{48}$ \\
\cline{4-5}
28 &           &         &        \{ 2, 4, 8, 15 \}              & ${\mathcal D}_{49}$ \\
\hline

29 & \multirow{2}*{$A_7$} & \multirow{2}*{$S_5$}
    &      \{ 7, 8, 9, 14, 15, 17, 18, 19, 20 \}            & ${\mathcal D}_{50}$ \\
\cline{4-5}
30 &           &         &     \{ 7, 8, 14, 18, 21 \}                 & ${\mathcal D}_{51}$ \\
\hline

31 & $A_7$   & $(A_4\times S_3)$:2
  &  \{ 2, 3, 5, 7, 8, 9, 10, 13, 14, 20, 23, 24, 28, 29, 30, 31, 32, 34 \}   & ${\mathcal D}_{52}$ \\
\hline

32 & \multirow{5}*{$A_8$} & \multirow{5}*{$2^3$:$L_3(2)$}
    &     \{ 1, 3, 7, 10, 15 \}     & ${\mathcal D}_{53}$ \\
\cline{4-5}
33 &           &         &       \{ 2, 5, 6, 7, 10, 12 \}              & ${\mathcal D}_{54}$ \\
\cline{4-5}
34 &           &         &       \{ 2, 3, 4, 5, 6, 7, 9, 11, 13, 14 \}              & ${\mathcal D}_{55}$ \\
\cline{4-5}
35 &           &         &       \{ 1, 3, 4, 8, 9, 11, 13, 14, 15 \}              & ${\mathcal D}_{56}$ \\
\cline{4-5}
36 &           &         &         \{ 1, 8, 12, 14 \}            & ${\mathcal D}_{57}$ \\
\hline

37 & \multirow{2}*{$A_8$} & \multirow{2}*{($A_5\times 3$):2}
    &      \{ 1, 2, 3, 10, 19, 23, 34, 37, 41, 44, 48, 52 \}                        & ${\mathcal D}_{58}$ \\
\cline{4-5}
38 &           &         &    \{ 1, 3, 19, 20, 25, 31, 34, 36, 41, 43, 44, 45 \}     & ${\mathcal D}_{59}$ \\
\hline

39 & \multirow{2}*{$S_8$} & \multirow{2}*{$S_5\times S_3$}
    &        \{ 1, 17, 18, 19, 21, 26, 39, 40, 41, 42, 45, 50 \}                   & ${\mathcal D}_{60}$ \\
\cline{4-5}
40 &           &         &   \{ 2, 6, 9, 10, 15, 16, 20, 22, 37, 47, 52, 53 \}      & ${\mathcal D}_{61}$ \\
\hline

41 & \multirow{4}*{$A_9$} & \multirow{4}*{$S_7$}
    &        \{ 6, 7, 10, 15, 19, 33 \}             & ${\mathcal D}_{62}$ \\
\cline{4-5}
42 &           &         &       \{ 2, 4, 6, 9, 13, 14, 20, 21, 22, 25, 28, 30, 31, 33, 35, 36 \}                      & ${\mathcal D}_{63}$ \\
\cline{4-5}
43 &           &         &          \{ 2, 4, 14, 20, 23, 35 \}                   & ${\mathcal D}_{64}$ \\
\cline{4-5}
44 &           &         &           \{ 5, 7, 20, 21, 24, 29 \}                  & ${\mathcal D}_{65}$ \\
\hline

45 & \multirow{3}*{$A_9$} & \multirow{3}*{$L_2(8)$:$ 3$}
    &    \{ 2, 10, 12, 16, 33, 34, 40, 45, 54, 65, 70, 71, 74, 91, 95, 102, 110, 118 \}
            & ${\mathcal D}_{66}$ \\
\cline{4-5}
46 &           &         &
 \{ 1, 2, 3, 17, 21, 22, 35, 39, 54, 66, 69, 79, 91, 97, 100, 103, 109, 113 \}
         & ${\mathcal D}_{67}$ \\
\cline{4-5}
47 &           &         &        \{ 1, 4, 23, 24, 30, 36, 37, 39, 53, 65, 67, 80, 82, 92, 95, 98, 109, 117 \}
        & ${\mathcal D}_{68}$ \\
\hline

\multirow{2}*{48} & \multirow{2}*{$S_9$} & \multirow{2}*{$3^3:(2\times S_4)$}
  &      \{ 1, 7, 24, 26, 31, 36, 60, 67, 88, 99, 102, 108, 118, 120, 122, 123, 127,
            & \multirow{2}*{${\mathcal D}_{69}$} \\

 &   &    &
     130, 156, 162, 171, 173, 178, 185, 196, 202, 208, 212, 220, 251, 267, 272 \}           &  \\
\hline

49 & \multirow{3}*{$A_{10}$} & \multirow{3}*{$S_8$} &  \{ 5, 8, 14, 15, 16, 17, 18, 26, 32, 34, 37, 40 \}    & ${\mathcal D}_{70}$ \\
\cline{4-5}
50 &           &         &         \{ 4, 7, 9, 12, 18, 25, 27, 28, 31, 33, 34, 40 \}               & ${\mathcal D}_{71}$ \\
\cline{4-5}
51 &           &         &        \{ 7, 12, 13, 15, 16, 20, 23, 34, 35, 39, 40, 41 \}                & ${\mathcal D}_{72}$ \\
\hline

52 & \multirow{2}*{$A_{10}$}    & \multirow{2}*{$(A_7\times 3)$:2}      &
\{ 3, 4, 7, 9, 17, 20, 22, 25, 26, 29, 45, 64, 67, 72, 83, 87, 93, 111 \}
             & ${\mathcal D}_{73}$ \\
\cline{4-5}
53 &           &         &     \{ 1, 11, 22, 23, 32, 35, 36, 47, 51, 53, 55, 59, 67, 71, 83, 88, 94, 110 \}                  & ${\mathcal D}_{74}$ \\
\hline

54 & \multirow{14}*{$A_{10}$} & \multirow{14}*{$(A_5\times A_5)$:4} &
           \{ 11, 16, 30, 56, 98, 112 \}            & ${\mathcal D}_{75}$ \\
\cline{4-5}
\multirow{2}*{55} &           &         &             \{ 3, 6, 12, 14, 18, 23, 24, 25, 26, 28, 30, 33, 34, 37, 38, 43, 46, 49,
    53,
     & \multirow{2}*{${\mathcal D}_{76}$} \\

 &           &         &         74, 81, 85, 93, 95, 96, 97, 98, 99, 100, 101, 106, 110, 114, 115, 121, 122 \}              & \\

\cline{4-5}
56 &           &         &           \{ 33, 56, 63, 67, 83, 111 \}                        & ${\mathcal D}_{77}$ \\
\cline{4-5}
57 &           &         &           \{ 24, 31, 68, 99, 113, 123 \}                        & ${\mathcal D}_{78}$ \\
\cline{4-5}
58 &           &         &        \{ 2, 16, 21, 25, 26, 31, 55, 58, 71, 73, 79, 83, 85, 93, 103, 114 \}                           & ${\mathcal D}_{79}$ \\
\cline{4-5}
59 &           &         &              \{ 22, 30, 31, 80, 97, 113 \}                     & ${\mathcal D}_{80}$ \\
\cline{4-5}
60 &           &         &            \{ 1, 7, 76, 83, 119, 122 \}                       & ${\mathcal D}_{81}$ \\
\cline{4-5}
61 &           &         &          \{ 7, 11, 12, 22, 34, 46, 51, 56, 62, 68, 82, 87, 93, 94, 107, 121 \}                         & ${\mathcal D}_{82}$ \\
\cline{4-5}
\multirow{2}*{62} &           &         &      \{ 1, 8, 13, 21, 23, 24, 28, 32, 35, 37, 43, 46, 47, 50, 52, 59, 61, 64, 65, 68,                             & \multirow{2}*{${\mathcal D}_{83}$} \\

 &           &         &     75, 76, 77, 78, 85, 86, 87, 95, 104, 105, 107, 113, 120, 123, 124, 125 \}       &  \\

\cline{4-5}
\multirow{2}*{63} &           &         &          \{ 4, 5, 6, 14, 16, 19, 21, 28, 36, 43, 53, 54, 58, 59, 60, 62, 63, 64, 66,
 & \multirow{2}*{${\mathcal D}_{84}$} \\

 &           &         &       68, 69, 71, 73, 74, 75, 94, 95, 97, 104, 107, 108, 109, 110, 120, 121, 123 \}       &  \\

\cline{4-5}
64 &           &         &            \{ 1, 17, 19, 26, 37, 43, 51, 55, 59, 60, 64, 68, 69, 87, 89, 121 \}                       & ${\mathcal D}_{85}$ \\
\cline{4-5}
\hline

65 & $S_{10}$  & $S_7\times S_3$ &    \{ 4, 5, 7, 22, 29, 37, 40, 42, 47, 64, 69, 75, 78, 80, 88, 92, 109, 115 \}                      & ${\mathcal D}_{86}$ \\
\hline

66 & $S_{10}$  & $(S_5\times S_5)$:2 &      \{ 1, 9, 23, 50, 112, 115 \}               & ${\mathcal D}_{87}$ \\
\hline

\end{tabular}
}
\end{center}
\end{table*}

\newpage

\textbf{Proof of Theory \ref{thm2}.}

In the subsection 4.1, we calculated all possible design 5-tuple parameters that satisfy the basic requirements.
With the help of group theory software $Magma$, we filtered out those parameters that cannot form a design.

When $T = A_n$($5 \leq n \leq 10$) and $G$ acts on $n$ points,
from the content of section 3,
we know that $G = A_n$ or $G_n$, and according to Theorem \ref{thm1},
there are 21 designs in total listed in Table \ref{tab:v=n}, all of which are full designs.
When $T = A_n$($5 \leq n \leq 10$) and $G$ acts on $v$ $(v\neq n)$ points,
all the designs listed in \ref{tab:v<>n}.
In the subsection 4.2, we listed the basic blocks of all designs up to isomorphism.

  This completes the proof of Theorem \ref{thm2}.
$\hfill\square$

\section{Conclusion and Future Work}

  This paper systematically investigates non-trivial flag-transitive and point-primitive $2$-$(v,k,\lambda)$
 designs with the alternating group $A_{n}$($5 \leq n \leq 10$) as the socle of their automorphism groups.

  This result resolves a key problem in the classification of designs under this group action, significantly enriches the 2-design taxonomy,
  and offers methodological references for classifying designs under other alternating group, sporadic or exceptional simple groups.

{\bf Funding}

This work is supported by the National Natural Science Foundation of China(Nos:11801092,12271173),
Guangdong Basic and Applied Basic Research Foundation(No:2025A1515012072).

\end{document}